\documentclass[11pt]{article}

\usepackage[utf8]{inputenc}
\usepackage{leftidx}
\usepackage{amsmath}
\usepackage{amscd}
\usepackage{amssymb}
\usepackage{rotating}
\usepackage{amsfonts}
\usepackage{amsthm}
\usepackage{verbatim}
\usepackage{stmaryrd}
\usepackage{mathrsfs}
\usepackage{color}

\usepackage[top=1in, bottom=1in, left=1.25in, right=1.25in]{geometry}

\usepackage{titlesec}
       \titleformat{\chapter}[display]
             {\normalfont\Large\bfseries}{\thechapter}{11pt}{\Large}
       \titleformat{\section}
             {\normalfont\large\bfseries}{\thesection}{11pt}{\large}
       \titlespacing*{\chapter}{0pt}{0pt}{15pt} 
       \titlespacing*{\section}{0pt}{3.5ex plus 1ex minus .2ex}{2.3ex plus .2ex}

    \usepackage[titletoc]{appendix}

\input xy
\xyoption{all}

\allowdisplaybreaks[1]

\newcommand{\pqed}{\hfill\qedsymbol\\}

\newcommand{\id}{\mathrm{id}}

\newcommand{\Spec}{\mathrm{Spec}}

\newtheorem{theorem}{Theorem}[section]
\newtheorem{theorem/definition}{Theorem/Definition}[section]

\newtheorem{proposition}[theorem]{Proposition}
\newtheorem{lemma}[theorem]{Lemma}

\newtheorem{corollary}[theorem]{Corollary}

\newtheorem{conjecture}[theorem]{Conjecture}

\theoremstyle{remark}
\newtheorem{remark}[theorem]{Remark}

\theoremstyle{definition}
 \newtheorem{example}[theorem]{Example}
\newtheorem{definition}[theorem]{Definition}

\newtheorem{question}[theorem]{Question}

\begin{document}

\title
{\large{\textbf{Towards a  specialization map  modulo semi-orthogonal decompositions}}}
\author{\normalsize Xiaowen Hu}
\date{}
\maketitle

\begin{abstract}
We propose a conjecture on the existence of a specialization map for derived categories of smooth proper varieties modulo semi-orthogonal decompositions, and verify it for K3 surfaces and abelian varieties.
\end{abstract}

\section{Introduction}
In this paper we study the following question: given a family of smooth projective varieties over, say, a punctured disc, and the knowledge of their bounded derived category of coherent sheaves, what can we say about the derived category of the limit fiber?

One motivation is the well-known conjecture of Dubrovin which  predicts that a smooth projective variety has semisimple quantum cohomology if and only if its derived category of coherent sheaves has a full exceptional collection (see \cite{Dub98}, and also \cite{Bay04}).  Since quantum cohomology is deformation invariant, it suggests that the property of having a full exceptional collection is also invariant under deformations. In \cite{Hu18} we showed that this is true locally; more precisely, given a  smooth proper scheme $\mathcal{X}$ over a locally noetherian scheme $S$, if for one fiber $\mathcal{X}_{s_0}$, $\mathrm{D}^{\mathrm{b}}(\mathcal{X}_{s_0})$ has a full exceptional collection, then so does the geometric fibers in an open neighborhood. It remains to investigate, with the additional hypothesis that  $S$ is connected, whether $\mathrm{D}^{\mathrm{b}}(\mathcal{X}_{s})$ has a full exceptional collection for each fiber $\mathcal{X}_s$. This reduces to the following :
\begin{question}\label{ques-1}
Let $R$ be a discrete valuation ring, $K$ its fraction field, and $k$ its residue field. Denote $S=\Spec(R)$,  the generic point of $S$ by $\eta$, and the closed point of $S$ by 0. Let $X$ be a smooth projective scheme over $S$. Suppose $\mathrm{D}^{\mathrm{b}}(X_{\eta})$ has a full exceptional collection. Then does $\mathrm{D}^{\mathrm{b}}(X_0)$ has a full exceptional collection?
\end{question}
  Now, given a field $k$, we consider the abelian group freely generated by the equivalence classes of derived categories of coherent sheaves of smooth projective varieties over $k$, and then modulo the relation of the form 
\begin{equation}\label{eq-semi-1}
[T]=[S_1]+...+[S_n]
\end{equation}
if there is a semi-orthogonal decomposition
\[
T=\langle S_1,...,S_n\rangle.
\]
We call the resulting group the \emph{Grothendieck group of strictly geometric triangulated categories over $k$}, and denote it by $K_0(\mathrm{sGT}_k)$. For brevity we denote the class of $\mathrm{D}^{\mathrm{b}}(X)$ in $K_0(\mathrm{sGT}_k)$ by $[X]$. If $\mathrm{D}^{\mathrm{b}}(X)$ has a full exceptional collection of length $n$, then
\[
[X]=n[\Spec(k)].
\]
 So a question weaker than \ref{ques-1} is, with the same  hypothesis, whether $[X_0]=n[\Spec(k)]$, where $n$ is the length of the full exceptional collection of $\mathrm{D}^{\mathrm{b}}(X_\eta)$. Furthermore, for this weaker question, one can weaken the hypothesis, i.e., instead of assuming $\mathrm{D}^{\mathrm{b}}(X_{\eta})$ has a full exceptional collection, we now only assume $[X_\eta]=n[\Spec(K)]$. More generally, we propose the following conjecture.

 \begin{conjecture}\label{conj-1}
There is a natural group homomorphism
\[
\rho_{\mathrm{sgt}}:
K_0(\mathrm{sGT}_K)\rightarrow K_0(\mathrm{sGT}_k).
\]
 \end{conjecture}
If such a map $\rho_{\mathrm{sgt}}$ does exist, we call it the specialization map of Grothendieck group of strictly geometric triangulated categories. The validity of conjecture \ref{conj-1} would be an evidence to a positive answer to question \ref{ques-1}. This conjecture is  inspired also by \cite{NS17} and \cite{KT17}, where the existence of certain specialization maps are used to show that stable rationaly and rationality are closed properties in a smooth proper family. For example, in \cite{NS17}, it is shown that there is a natural group homomorphism
\[
\rho_{\mathrm{Var}}:
K_0(\mathrm{Var}_K)\rightarrow K_0(\mathrm{Var}_k).
\]
It is not hard to see that there is a canonical surjective homomorphism (see section 2)
\[
\mu:
K_0(\mathrm{Var}_k)/(\mathbb{L}-1)
\rightarrow K_0(\mathrm{sGT}_k).
\]
It should be believed that $\mu$ is not an isomorphism, but this problem seems still open.

In this paper we propose a definition of the map $\rho_{\mathrm{sgt}}$, and verify the well-definedness for K3 surfaces and abelian varieties. 

A more natural object to study than $K_0(\mathrm{sGT}_K)$ is the group generated by the admissible subcategories of derived categories of coherent sheaves of smooth projective varieties, modulo the same kind of relations (\ref{eq-semi-1}), and one can propose a conjecture parallel to conjecture \ref{conj-1}. A closedly related notion is the Grothendieck ring of pre-triangulated categories introduced in \cite{BLL04}.\\

\noindent\textbf{Acknowledgement}\quad
I am grateful to Lei Zhang,  Zhan Li, Qingyuan Jiang, Ying Xie, Shizhuo Zhang,  Xin Fu, Feng Qu and Lei Song for helpful discussions. Part of this paper  is inspired by a workshop  in SUSTech organized by Zhan Li. 

This work is supported by 34000-31610265, NSFC 34000-41030364 and 34000-41030338.

\section{Definitions and the conjecture}
Throughout this section, denote by $k$  a field of characteristic zero, $R=k[[t]]$, and denote by $K$ the fraction field of $R$.  For a smooth  proper scheme $X$ over  $K$, an \emph{snc model}  of $X$ is a proper scheme $\mathcal{X}$  over $R$ with the properties that $\mathcal{X}$ is regular, $\mathcal{X}_K$ is isomorphic to $X$ and the special fiber $\mathcal{X}_0=\mathcal{X}\times_R k$ is an snc divisor of $\mathcal{X}$, which means
\[
\mathcal{X}_0=\bigcup_{i=1}^{n}m_i D_i
\]
 as divisors, where $D_i$ is an irreducible smooth proper scheme over $k$, $m_i$ are positive integers, and if one writes $D_I=\bigcap_{i\in I}D_i$ for $I\subset \{1,...,n\}$, then  $\dim D_I=\dim X+1-|I|$ for all subsets $I$ of $\{1,...,n\}$. Our only use of the assumption of the characteristic 0 is that such fields admit resolution of singularities and the weak factorization theorems hold in this case (\cite{AKMW02}, \cite{Wlo03}, \cite{AT16}). In particular, snc models always exists.
\begin{definition}
A $k$-linear triangulated category $T$ is \emph{geometric}, if there is a smooth proper scheme $Y$ over $k$ such that $T$ is equivalent to an admissible triangulated subcategory of $\mathrm{D}^{\mathrm{b}}(Y)$.  A $k$-linear triangulated category $T$ is \emph{strictly geometric}, if there is a smooth proper scheme $Y$ over $k$ such that $T$ is equivalent to  $\mathrm{D}^{\mathrm{b}}(Y)$.
\end{definition}

Let $K_0(\mathrm{GT}_k)$ (resp., $K_0(\mathrm{sGT}_k)$) be the quotient of the free abelian group generated by the equivalence classes of geometric (resp., strictly geometric) triangulated categories modulo the relations of the form
\[
[T]=[T_1]+...+[T_n]
\]
if there is a semi-orthogonal decomposition 
\[
\langle S_1,...,S_n\rangle
\]
 of $T$ such that  $S_i$ is equivalent to $T_i$ for $1\leq i\leq n$.  In particular, the class of the  zero category is equal to zero.
 
Recall the following two theorems of Orlov  on the semi-orthogonal decomposition of projective bundles and blow-ups (see \cite{Orl92}, or \cite[chapter 8]{Huy06}).
\begin{theorem}\label{so-pb}
Let $Y$ be a smooth projective variety over $k$, $E$ be a vector bundle of rank $r$ over $Y$, and $\pi: \mathbb{P}(E)\rightarrow Y$ the projective bundle. Then there is a semi-orthogonal decomposition
\begin{equation}\label{so-pb-1}
\mathrm{D}^{\mathrm{b}}(\mathbb{P}(E))=\langle \pi^* \mathrm{D}^\mathrm{b}(Y)\otimes \mathcal{O}(a),...,\pi^* \mathrm{D}^\mathrm{b}(Y)\otimes \mathcal{O}(a+r-1)\rangle
\end{equation}
for every integer $a$.
\end{theorem}

\begin{theorem}\label{so-blowup}
Let $X$ be a smooth projective variety over $k$ and $Y$ a smooth closed subvariety of $X$ of codimension $c\geq 2$, $\mathrm{Bl}_Y X$ the blowup of $X$ along $Y$. Then there is a semi-orthogonal decomposition
\begin{equation}\label{so-blowup-1}
\mathrm{D}^{\mathrm{b}}(\mathrm{Bl}_Y X)=\langle \mathcal{D}_{-c+1},...,\mathcal{D}_{-1},\mathrm{D}^{\mathrm{b}}(X)\rangle
\end{equation}
such that $\mathcal{D}_i$ is equivalent to $\mathrm{D}^b(Y)$ for $-c+1\leq i\leq -1$.
\end{theorem}

Denote by $K_0(\mathrm{Var}_k)$ the Grothendieck group of varieties over $k$.  Recall that $K_0(\mathrm{Var}_k)$ is the group generated by the isomorphism classes of smooth schemes over $k$ modulo the relations
\[
 [X]=[Y]+[U]
 \] 
 where $X$ is a smooth scheme over $k$, $Y$ is a closed subscheme of $X$ which is also smooth over $k$, and $U=X-Y$. The following theorem of Bittner \cite[theorem 3.1]{Bit04} gives an equivalent definition.
\begin{theorem}\label{thm-bit}
 Let $k$ be a field of characteristic zero. Then $K_0(\mathrm{Var}_k)$ is isomorphic to the
  group generated by the isomorphism classes of smooth proper schemes over $k$ modulo the relations
\[
 [X]-[Y]=[\mathrm{Bl}_Y X]-[E]
 \] 
 where $X$ is a smooth proper scheme over $k$, $Y$ is a smooth closed subscheme of $X$, $\mathrm{Bl}_Y X$ is the blow-up of $X$ along $Y$, and $E$ is the corresponding exceptional divisor on $\mathrm{Bl}_Y X$. 
\end{theorem}

\begin{corollary}\label{cor-mu}
Suppose $k$ is a field of characteristic zero. Then there is a natural surjective homomorphism of groups
\begin{equation}\label{eq-mu}
\mu_k: K_0(\mathrm{Var}_k)/[\mathbb{L}-1]\rightarrow K_0(\mathrm{sGT}_k)
\end{equation}
such that $\mu_k([X])=[\mathrm{D}^{\mathrm{b}}(X)]$.
\end{corollary}
Proof :  By theorem \ref{thm-bit}, it suffices to show 
\begin{equation}\label{eq-cor-mu-1}
[\mathrm{D}^\mathrm{b}(X)]-[\mathrm{D}^\mathrm{b}(Y)]=[\mathrm{D}^\mathrm{b}(\mathrm{Bl}_Y X)]-[\mathrm{D}^\mathrm{b}(E)]
\end{equation}
and
\begin{equation}\label{eq-cor-mu-2}
[\mathrm{D}^\mathrm{b}(\mathbb{P}^1)]=2[\mathrm{D}^\mathrm{b}(\Spec(k))].
\end{equation}
By (\ref{so-pb}), $[\mathrm{D}^\mathrm{b}(\mathbb{P}^n)]=(n+1)[\mathrm{D}^\mathrm{b}(\Spec(k))]$, thus (\ref{eq-cor-mu-2}) holds. Suppose the codimension of $Y$ in $X$ is $c$, then by (\ref{so-blowup}), 
\[
[\mathrm{D}^\mathrm{b}(\mathrm{Bl}_Y X)]=(c-1)[\mathrm{D}^\mathrm{b}(Y)]+[\mathrm{D}^\mathrm{b}(X)],
\]
and by (\ref{so-pb}), 
\[
[\mathrm{D}^\mathrm{b}(E)]=c[\mathrm{D}^\mathrm{b}(Y)],
\]
so (\ref{eq-cor-mu-1}) follows. Since $K_0(\mathrm{Var}_k)$ and $K_0(\mathrm{sGT}_k)$ both are generated by the isomorphism classes of smooth proper schemes over $k$, $\mu_k$ is surjective. \pqed

\begin{remark}
We have ignored the ring structure of $K_0(\mathrm{Var}_k)$. To obtain a ring structure on something like $K_0(\mathrm{sGT}_k)$ or $K_0(\mathrm{GT}_k)$, one need take into account the DG structrues (see \cite{BLL04}), and there is then a map like (\ref{eq-mu}).
\end{remark}

Now let $k$ and $K$ be the fields as defined at the beginning of this section. The following theorem is  \cite[prop. 3.2.1]{NS17}.

\begin{theorem}\label{thm-NS}
There is a unique group homomorphism
\[
\rho_{\mathrm{var}}: K_0(\mathrm{Var}_K)\rightarrow K_0(\mathrm{Var}_k)
\]
such that for a smooth proper scheme $X$ over $K$,  an snc model $\mathcal{X}$ of $X$ over $R$ with
\[
\mathcal{X}_k=\sum_{i\in I}n_i D_i,
\]
one has
\begin{equation}\label{eq-def-rho-var}
\rho_{\mathrm{var}}([X])=\sum_{\emptyset\neq J\subset I}(1- \mathbb{L} )^{|J|-1}[D_J^{\circ}],
\end{equation}
where $D_J=\bigcap_{j\in J}D_j$, and $D_J^{\circ}=D_J\backslash(\bigcup_{i\in I\backslash J}D_i)$.
\end{theorem}
\pqed

The homomorphism $\rho_{\mathrm{var}}$ is called the \emph{specialization map} of the Grothendieck group of varieties.

\begin{conjecture}\label{conj-spc-1}
 There are  natural  maps
 \[
 \rho_{\mathrm{gt}}: K_0(\mathrm{GT}_K)\rightarrow K_0(\mathrm{GT}_k)
 \]
 and
 \[
 \rho_{\mathrm{sgt}}: K_0(\mathrm{sGT}_K)\rightarrow K_0(\mathrm{sGT}_k).
 \]
 \end{conjecture}\pqed

In view of theorem \ref{thm-NS} and corollary \ref{cor-mu}, the conjecture for $K_0(\mathrm{sGT})$ means that there is a homomorphism $\rho_{\mathrm{sgt}}$ making the following diagram commutative
 \begin{equation}\label{graph-var-sgt}
 \xymatrix{
    K_0(\mathrm{Var}_K) \ar[r]^{\rho_{\mathrm{var}}} \ar[d]_{\mu_K} & K_0(\mathrm{Var}_k)  \ar[d]^{\mu_k} \\
    K_0(\mathrm{sGT}_K) \ar@{.>}[r]^{\rho_{\mathrm{sgt}}}  & K_0(\mathrm{sGT}_k),
 }
 \end{equation}
and since $\mu_K$ is surjective, such $\rho_{\mathrm{sgt}}$ is unique if it exists.

For a field $L$, denote by $\mathrm{M}_L$ the abelian group freely generated by the isomorphism classes of smooth proper schemes over $L$. 
 Set
\[
\mathbb{P}_{D_J}=\mathbb{P}(\mathcal{N}_{D_J/\mathcal{X}}).
\]
In particular, $\mathbb{P}_{D_i}=D_i$. We define a map
\[
\rho: \mathrm{M}_K\rightarrow K_0(\mathrm{sGT}_k)
\]
by
\[
\rho([X])=\sum_{\emptyset \neq J\subset I}(-1)^{|J|-1}[\mathrm{D}^{\mathrm{b}}(\mathbb{P}_{D_J})],
\]
or equivalently, by theorem \ref{so-pb},
\begin{equation}\label{eq-def-rho-sgt-1}
\rho([X])=\sum_{\emptyset \neq J\subset I}(-1)^{|J|-1}|J|\cdot[\mathrm{D}^{\mathrm{b}}(D_J)].
\end{equation}
By (\ref{eq-def-rho-var}) a simple computation shows that 
\[
\rho_{\mathrm{var}}([X])=\sum_{\emptyset \neq J\subset I}(-1)^{|J|-1}[\mathbb{P}_{D_J}].
\]
Therefore $\rho$ is a natural candidate for $\rho_{\mathrm{sgt}}$. In other words, conjecture \ref{conj-spc-1} for $\rho_{\mathrm{sgt}}$ reduces to the following.
\begin{conjecture}\label{conj-spc-2}
The homomorphism $\rho:\mathrm{M}_K\rightarrow K_0(\mathrm{sGT}_k)$ factors through the canonical surjective homomorphism $\mathrm{M}_K\twoheadrightarrow K_0(\mathrm{sGT}_K)$:
\[
\xymatrix{
  \mathrm{M}_K \ar[rr]^{\rho} \ar[rd] && K_0(\mathrm{sGT}_k)\\
  & K_0(\mathrm{sGT}_K) \ar@{.>}[ru] &.
}
\]
\end{conjecture}

  To prove the conjecture, one need to show:
\begin{enumerate}
  \item[(i)] given $X$ as a representative of its class $[\mathrm{D}^{\mathrm{b}}(X)]$ in $K_0(\mathrm{sGT}_K)$, $\rho([X])$ is independent of the choice of the snc model $\mathcal{X}$;
  \item[(ii)] $\rho([X])$ is independent of the choice of the representative $X$.
\end{enumerate}

In fact, (i) is needed for the well-definedness of $\rho$. We state it as follows.

\begin{theorem}\label{thm-indep-1}
$\rho([X])$  does not depend on the choice of $\mathcal{X}$.
\end{theorem}
Proof : One can show this by using the weak factorization theorem \cite{AKMW02}, \cite{Wlo03} and \cite{AT16}. The quickest way is to apply theorem \ref{thm-NS} and corollary \ref{cor-mu}.\pqed

I have no idea how to do step (ii) at present. In this paper I  only provide some evidence for it. More precisely, for some examples of derived equivalent smooth proper $K$-schemes $X$ and $X'$, I am going to verify
\begin{equation}\label{eq-rho-check-0}
\rho([X])=\rho([X']).
\end{equation}

The first kind of examples are birational derived equivalent $X$ and $X'$. 

\begin{lemma}\label{lem-property-rho}
Let $X$ be a smooth proper scheme over $K$.
\begin{enumerate}
  \item[(i)]  Let $Y$ be a smooth closed subscheme of $X$. Denote by $E$ the exceptional divisor on the blowup $\mathrm{Bl}_Y X$. Then $\rho([\mathrm{Bl}_{Y}X])=\rho([X])-\rho([Y])+\rho([E])$.
  \item[(ii)] Let $E$ be a vector bundle over $X$ of rank $r$. Then $\rho(\mathbb{P}(E))=r\rho([X])$.
\end{enumerate}
\end{lemma}
Proof : Use corollary \ref{cor-mu} and theorem \ref{thm-NS}. \pqed

\begin{example}[Standard flips]
Let $X$ be a smooth projective scheme over $K$ and $Y$ a smooth closed subscheme of $X$ of codimension $l+1$, such that $Y\cong \mathbb{P}^m$ and the normal bundle $N_{Y/X}\cong \mathcal{O}(-1)^{l+1}$. Then one can perform the standard flip and obtain a smooth projective scheme $X'$. By \cite[theorem 3.6]{BO95}, $X$ and $X'$ are derived equivalent.  By lemma \ref{lem-property-rho} one deduces that   
\[
\rho([X])+l \rho([\mathbb{P}^m])=\rho([\mathrm{Bl}_Y X])=\rho([X'])+l \rho([\mathbb{P}^m]),
\]
so
\[
\rho([X])=\rho([X']).
\]
\end{example}
\pqed

Similarly, one can also try to check (\ref{eq-rho-check-0}) for Mukai flops (\cite{Kaw02},\cite{Nam03}), and two non-isomorphic crepant resolutions of a Calabi-Yau 3-fold. In the following sections I will verify (\ref{eq-rho-check-0}) for K3 surfaces and abelian varieties, under some additional assumptions.

\section{Specialization map   K3 surfaces}
In this section we verify (\ref{eq-rho-check-0}) for derived equivalent $K3$ surfaces which have semistable degenerations over $R$. Throughout this section we consider only algebraic K3 surfaces.
\subsection{Mukai pairings and period mappings}
In this subsection we recall the Mukai pairing on K3 surfaces and its relation to derived equivalences (see e.g., \cite[chapter 4]{BBR09} and \cite[chapter 10]{Huy06}), and then introduce a corresponding notion of  period mapping.  

Let $X$ be a K3 surface over $\mathbb{C}$. The \emph{Mukai pairing} on $H^*(X,\mathbb{Z})$ is defined by
\[
\langle (\alpha_0,\alpha_1,\alpha_2),(\beta_0,\beta_1,\beta_2)\rangle:=a_1.\beta_1- \alpha_0.\beta_2- \alpha_2.\beta_0\in \mathbb{Z},
\]
where $\alpha_i,\beta_i\in H^{2i}(X,\mathbb{Z})$. The corresponding lattice is
\[
E_8(-1)^{\oplus 2}\oplus U^{\oplus 4}.
\]
Set
\begin{eqnarray*}
\widetilde{H}^{2,0}(X)=H^{2,0}(X),\ \widetilde{H}^{0,2}(X)=H^{0,2}(X),\\
\widetilde{H}^{1,1}(X)=H^{0}(X)\oplus H^{4}(X)\oplus H^{1,1}(X).
\end{eqnarray*}
The resulting weight two Hodge structure $\{H^{\mathrm{even}}(X,\mathbb{Z}),\widetilde{H}^{p,q}(X)\}$ is denoted by $\widetilde{H}(X,\mathbb{Z})$.

The following characterization of derived equivalent K3 surfaces is due to \cite{Muk87}, \cite{Orl97}. See also \cite[corollary 10.7, proposition 10.10]{Huy06}.
\begin{theorem}\label{thm-equiv-K3}
Two algebraic K3 surfaces $X$ and $Y$ over $\mathbb{C}$ are derived equivalent if and only if there is a Hodge isometry between $\widetilde{H}(X,\mathbb{Z})$ and $\widetilde{H}(Y,\mathbb{Z})$ with respect to the Mukai pairing. If $\Phi_P:\mathrm{D}^{\mathrm{b}}(X)\rightarrow \mathrm{D}^\mathrm{b}(Y)$ is an equivalence with kernel $P\in \mathrm{D}^{\mathrm{b}}(X\times Y)$, the induced map 
\[
\Phi_P^H: \widetilde{H}(X,\mathbb{Z})\rightarrow \widetilde{H}(Y,\mathbb{Z}),\ \alpha\mapsto q_* (\mathrm{ch}(P)\mathrm{td}(X\times Y)\cdot p^*\alpha)
\]
is a Hodge isometry, where $p:X\times Y\rightarrow X$, $q:X\times Y\rightarrow Y$ are the two projections.
\end{theorem}

As an analogue of the usual period domains, we introduce a notion to study the variation of $\widetilde{H}(X,\mathbb{Z})$.
\begin{definition}
Let $\mathbf{M}$ be the Mukai lattice $E_8(-1)^{\oplus 2}\oplus U^{\oplus 4}$, $Q(\cdot ,\cdot)$ the corresponding symmetric bilinear pairing on $H_{\mathbb{C}}$. The \emph{Mukai period domain} $D_M$ is defined to be the classifying space of the following data:
\begin{enumerate}
   \item[(i)] a filtration of complex subspaces $0=F^3\subset F^2\subset F^1\subset F^0=\mathbf{M}_{\mathbb{C}}$ of $\mathbf{M}_{\mathbb{C}}$, such that $\dim_{\mathbb{C}}(F^2)=1$, $\dim_{\mathbb{C}}(F^1)=23$;
   \item[(ii)] $Q(F^p,F^{3-p})=0$ for all $p$;
   \item[(iii)] $Q(v,\bar{v})>0$ for $v\in F^2$.
 \end{enumerate} 
\end{definition}

Notice that the  condition (iii) together with condition (ii) implies that $F^1\cap \overline{F^2}=0$, thus
induces a weight two integral Hodge structure on $\mathbf{M}$.

\begin{proposition}\label{prop-domain}
\begin{enumerate}
  \item[(i)] $D_M$ is an open subset (in the analytic topology) of a subvariety of a flag variety;
  \item[(ii)] For a family of K3 surface $\mathcal{X}\rightarrow S$, where $S$ is a simply connected complex manifold, and an isomorphism $\widetilde{H}^*(\mathcal{X}_0,\mathbb{Z})\cong \mathbf{M}$ as lattices for some point $0$ of $S$, there is a canonical holomorphic  map $\phi:S\rightarrow D_M$, 
  such that 
  \[
  \widetilde{H}(\mathcal{X}_s,\mathbb{Z})\cong \phi(s), 
  \]
  for any  point $s$ of $S$.
\end{enumerate}
\end{proposition}
Proof: Both statements follow from the usual argument for the period map of unpolarized K3 surfaces, see \cite[chapter 6]{Huy16}. For example, $\coprod_{s\in S}H^0(\mathcal{X}_s, \Omega_{\mathcal{X}_s}^2)=\mathrm{R}^0\pi_*\Omega_{\mathcal{X}/S}^2$  is a holomorphic subbundle of $\bigoplus_{i=0}^{4}\mathrm{R}^i\pi_*\mathbb{Z}\otimes_{\mathbb{Z}}\mathcal{O}_S$, so $\phi$ is holomorphic.
\pqed

More generally, for 
a family of K3 surface $\mathcal{X}\rightarrow S$, where $S$ is a complex manifold which is not necessarily simply connected, and an isomorphism $\widetilde{H}^*(\mathcal{X}_0,\mathbb{Z})\cong \mathbf{M}$ as lattices for some $0\in S$, there is a canonical holomorphic  map $\phi:S\rightarrow \Gamma\backslash D_M$, where $\Gamma=\mathrm{Aut}_{\mathbb{Z}}(\mathbf{M},Q)$, the group of automorphisms of the lattice $(\mathbf{M},Q)$, or even the image of $\pi_1(S)$ in $\mathrm{Aut}_{\mathbb{Z}}(\mathbf{M},Q)$. However, the quotient  $\Gamma\backslash D_M$ is not Hausdorff, as remarked in \cite[p. 104]{Huy16}.

\subsection{Degeneration of K3 surfaces}
We first recall the theorem on the degeneration of K3 surfaces due to Kulikov \cite{Kul77} (see also \cite{PP81},  \cite{Fri84}).
\begin{theorem}\label{thm-deg-K3}
Let $\pi:\mathcal{X}\rightarrow \Delta$ be a semistable degeneration of K3 surfaces. Then there exists a birational modification of this semistable degeneration, such that the restriction of $\pi$ to $\Delta^*=\Delta\backslash\{0\}$ remains unchanged, and $K_{\mathcal{X}}$ becomes trivial. After such a modification, the degenerate fiber $\pi^{-1}(0)=\mathcal{X}_0$ can be one of the following types:
\begin{enumerate}
  \item[(I)] $\mathcal{X}_0$ is a smooth K3 surface;
  \item[(II)] $\mathcal{X}_0=\bigcup_{i=0}^{r}V_i$, $r\geq 1$, $V_0$ and $V_r$ are rational surface, $V_i$ are ruled elliptic surfaces for $1\leq i\leq r-1$, $V_{i}\cap V_{j}=\emptyset$ for $|i-j|>1$, and  $V_{i}\cap V_{j}$ is an elliptic curve for  $|j-i|=1$ and is a section of the ruling on $V_i$, if $V_i$ is a ruled elliptic surface;
  \item[(III)] $\mathcal{X}_0=\bigcup V_i$, and each $V_i$ is a smooth rational surface, with all the double curves rational, and the dual graph is a triangulation of $S^2$.
\end{enumerate}
Moreover, the three types of degenerations are characterized by the monodromy action $T$ on $H^2(X_t,\mathbb{Z})$, $0\neq t\in \Delta$:
\begin{enumerate}
  \item[(I)] $T=\id$;
  \item[(II)] $T-\id\neq 0$, $(T-\id)^2=0$;
  \item[(III)] $(T-\id)^2\neq 0$, $(T-\id)^3=0$.
\end{enumerate}
\end{theorem}

In the following we say that a semistable degeneration of K3 surfaces with $K_{\mathcal{X}}$ trivial, is of type (I), (II) or (III), if it is of the corresponding type described above. 

\begin{proposition}\label{prop-deg-type-ii}
Let $\pi:\mathcal{X}\rightarrow \Delta$ be a type (II) semistable degeneration of K3 surfaces. Denote by $LH^i(\mathcal{X}_0)$ the limit Hodge structure on $H^i(\mathcal{X}_t)$. Denote by $E$ the elliptic curve which is isomorphic to the base elliptic curves of the ruled elliptic  surfaces appearing in $\mathcal{X}_0$. Then 
\begin{enumerate}
  \item[(i)] $W_1 H^2(\mathcal{X}_0)\cong H^1(E)^{\oplus r}$ as integral pure Hodge structures;
  \item[(ii)] $W_1 H^2(\mathcal{X}_0)\cong W_1 LH^2(\mathcal{X}_0)$ as integral pure Hodge structures.
  
\end{enumerate}
\end{proposition}
Proof: Notice that in general the pure graded pieces in a mixed Hodge structure are \emph{rational} Hodge structures. However in our case there are natural integral Hodge structures on $W_1 H^2(\mathcal{X}_0)$ and $W_1 LH^2(\mathcal{X}_0)$ inducing the rational ones as we will see. So it suffices to show the isomorphisms as rational Hodge structures.

(i) Let $\mathcal{X}_0=\bigcup_{i=0}^r V_i$.  For $j=1,...,r$, let $U_j=\bigcup_{i=0}^{j-1}V_i$ and $U'_j=\bigcup_{i=j}^{r}V_i$, and $D_j=V_{j-1}\cap V_j$. The Mayer-Vietoris exact sequence 
\[
...\rightarrow H^{k-1}(D_j)\rightarrow H^k(\mathcal{X}_0)\rightarrow H^k(U_j)\oplus H^k(U'_j)\rightarrow H^k(D_j)\rightarrow...
\]
provides an extension of pure Hodge structures
\[
0\rightarrow H^1(D_j)\rightarrow W_1 H^2(\mathcal{X}_0)\rightarrow W_1 H^2(U_j)\oplus W_1 H^2(U'_j)\rightarrow 0.
\]
By induction on $r$, $W_1 H^2(\mathcal{X}_0)$ is a successive extension of $H^1(D_1),...,H^1(D_r)$. But we can choose different orders of the cuts of $\mathcal{X}_0$, which give the splittings. Hence there is a canonical isomorphism $W_1 H^2(\mathcal{X}_0)\cong \bigoplus_{i=1}^r H^1(D_i)$ of Hodge structures.

(ii) By \cite[lemma 3.6]{Fri84}, the Clemens-Schmidt sequence
\[
H_4(\mathcal{X}_0)\rightarrow H^2(\mathcal{X}_0)\rightarrow LH^2(\mathcal{X}_0)\xrightarrow{N=T-1} LH^2(\mathcal{X}_0)
\]
is exact over $\mathbb{Z}$, and is an exact sequence of mixed Hodge structures. Since $N(W_1L H^2(\mathcal{X}_0))=0$, we have $W_1 H^2(\mathcal{X}_0)\cong W_1 LH^2(\mathcal{X}_0)$ as Hodge structures. \pqed

\subsection{The specialization map}
\begin{proposition}\label{prop-openness}
Let $R$ be an integral domain, $K$ the fraction field of $R$. Let
$\mathcal{X}$ and $\mathcal{Y}$ be smooth projective schemes over $R$, and 
 $X=\mathcal{X}_K$, $Y=\mathcal{Y}_K$. 
Suppose  $\Phi:\mathrm{D}^{\mathrm{b}}(X)\rightarrow\mathrm{D}^{\mathrm{b}}(Y)$ is an exact functor which is an equivalence of triangulated categories. 
  Then there exists $0\neq r\in R$  and $\mathcal{P}\in \mathrm{D}^{\mathrm{b}}(\mathcal{X} \times_R \mathcal{Y}) $ such that  for every point $s\in \Spec(R[\frac{1}{r}])$, the Fourier-Mukai transform 
\[
\Phi_{\mathcal{P}_s}: \mathrm{D}^{\mathrm{b}}(\mathcal{X}_s)\rightarrow \mathrm{D}^{\mathrm{b}}(\mathcal{Y}_s) 
\]
induced by $\mathcal{P}_s$ is an equivalence, wherer $\mathcal{X}_s$ and $\mathcal{Y}_s$ are the fiber over the point $\iota_s: \Spec(\kappa(s))\rightarrow \Spec(R[\frac{1}{r}])$, and $\mathcal{P}_s=\mathbf{L}\iota_s^* \mathcal{P}$, and moreover, $\Phi_{\mathcal{P}_s}=\Phi$.
\end{proposition}
Proof: By \cite[theorem 3.2.2]{Orl03}, there exists  $P\in \mathrm{D}^{\mathrm{b}}(X\times_K Y)$ such that $\Phi=\Phi_{P}$.
Shrinking $\Spec(R)$ if necessary, one can find $\mathcal{P}\in \mathrm{D}^{\mathrm{b}}(\mathcal{X} \times_R \mathcal{Y})$ such that  $\mathcal{P}_K=P$. Let $L$ be a very ample line bundle over $X$. Shrinking $\Spec(R)$ if necessary, there exists a relatively ample line bundle $\mathcal{L}$ on $\mathcal{X}$ over $R$, such that $\mathcal{L}$ restricts to $L$. Set 
\[
\mathcal{E}=\bigoplus_{i=0}^{d}\mathcal{L}^{\otimes i},
\]
where $d=\dim X=\dim Y$. By \cite[theorem 4]{Orl09}, $\mathcal{E}_s$  is a \emph{classical generator} of $\mathrm{D}^{\mathrm{b}}(\mathcal{X}_s)$, namely, the smallest triangulated subcategory of $\mathrm{D}^{\mathrm{b}}(\mathcal{X}_s)$ containing $\mathcal{E}_s$ and closed under isomorphisms and taking direct summands, is $\mathrm{D}^{\mathrm{b}}(\mathcal{X}_s)$.

 Set
\[
Q=P^{\vee}\otimes^{\mathbf{L}}p_2^* \omega_{Y}[d],\
\mathcal{Q}=\mathcal{P}^{\vee}\otimes^{\mathbf{L}}p_2^* \omega_{\mathcal{Y}/R}[d].
\]
Then $Q=\mathcal{Q}_K$, and for every point $s \in \Spec(R)$, the Fourier-Mukai transform $\Phi_{\mathcal{Q}_s}:\mathrm{D}^{\mathrm{b}}(\mathcal{Y}_s)\rightarrow \mathrm{D}^{\mathrm{b}}(\mathcal{X}_s)$ is a left adjoint of $\Phi_{\mathcal{P}_s}:\mathrm{D}^{\mathrm{b}}(\mathcal{X}_s)\rightarrow \mathrm{D}^{\mathrm{b}}(\mathcal{Y}_s)$.
Moreover, by hypothesis, $\Phi_{\mathcal{Q}_K}=\Phi_{Q}$ is an inverse of $\Phi_{\mathcal{P}_K}=\Phi_{P}$. Thus the adjoint map
\[
\Phi_Q\circ \Phi_P(\mathcal{E}_K)\rightarrow \mathcal{E}_K
\]
is an isomorphism. By the semi-continuity theorem (for perfect complexes, \cite[7.7.5]{EGAIII}),  shrinking $\Spec(R)$ if necessary, for every point $s \in \Spec(R)$,  the adjoint map
\[
\Phi_{\mathcal{Q}_s}\circ \Phi_{\mathcal{P}_s}(\mathcal{E}_s)\rightarrow \mathcal{E}_s
\]
is an isomorphism. By induction on the generating time of the objects of $\mathrm{D}^{\mathrm{b}}(\mathcal{X}_s)$ with respect ot $\mathcal{E}_s$, this implies that the adjoint morphism of functors 
\[
\Phi_{\mathcal{Q}_s}\circ \Phi_{\mathcal{P}_s}\rightarrow \id_{\mathrm{D}^{\mathrm{b}}(\mathcal{X}_s)}
\]
is an isomorphism. Thus $\Phi_{\mathcal{P}_s}$ is fully faithful. Finally starting from a very ample line bundle on $Y$, and shrinking $\Spec(R)$ if necessary, we find that $\Phi_{\mathcal{Q}_s}$ is also fully faithful, and we are done.\pqed

\begin{theorem}
Let $R$ be the local ring $\mathbb{C}[T]_{(T)}$, and $K=\mathbb{C}(T)$. Let $X$ and $Y$ be smooth projective surfaces  over $K$ with trivial canonical bundles. Suppose that $X$ and $Y$ are derived equivalent, and both have semistable degenerations over $R$. Then $\rho([X])=\rho([Y])$ in $K_0(\mathrm{sGT}_\mathbb{C})$.
\end{theorem}
Proof: Let $\mathcal{X}_R$ and $\mathcal{Y}_R$ be semistable degenerations of $X$ and $Y$ over $R$, respectively. Denote the point $(T)$ of $\Spec(R)$ and the point $(T)$ of $\Spec(\mathbb{C}[T])$, both by 0.  Then there exists an affine open subset $U$ of $\Spec(\mathbb{C}[T])$ and schemes $\mathcal{X}$ and $\mathcal{Y}$ over $U$, such that
\begin{enumerate}
  \item[(i)] restricting to $U\backslash\{0\}$, $\mathcal{X}$ and $\mathcal{Y}$  are smooth, and each geometric fiber is a K3 surface;
  \item[(ii)] the base changes of $\mathcal{X}$ and $\mathcal{Y}$ to $\Spec(R)$ are isomorphic to $\mathcal{X}_R$ and $\mathcal{Y}_R$, respectively.
\end{enumerate}
 By proposition \ref{prop-openness}, there is an open subset $V$ of $U$ containing $0$, and 
\[
\mathcal{P}\in \mathrm{D}^{\mathrm{b}}(\mathcal{X}\times_U \mathcal{Y}\times_U (V-\{0\}))
\]
  such that the  Fourier-Mukai transform 
 \[
 \Phi_{\mathcal{P}_t}: \mathrm{D}^{\mathrm{b}}(\mathcal{X}_t)\rightarrow \mathrm{D}^{\mathrm{b}}(\mathcal{Y}_t)
 \]
 is an equivalence, for all $t\in V-\{0\}$. Without loss of generality we assume $V=U$. Consider the analytic topology of $U$. Taking an open disk $\Delta$ of $U$ containing 0, and consider $\mathcal{X}$ and $\mathcal{Y}$ restricting over $\Delta$, we can apply the result of the previous subsections to study the fiber $\mathcal{X}_0$ and $\mathcal{Y}_0$.
By theorem \ref{thm-indep-1}, birational modifications preserving $\mathcal{X}-\mathcal{X}_0$ does not change $\rho([X])$. So by the first statement of theorem \ref{thm-deg-K3}, we can assume $K_{\mathcal{X}_\Delta}$ and $K_{\mathcal{Y}_\Delta}$ trivial, such that $\mathcal{X}_0$ and $\mathcal{Y}_0$ are described by theorem loc. cit.

For a point $t\in \Delta-0$, let $\Phi^H_{\mathcal{P}_t}: H^*(\mathcal{X}_t)\rightarrow H^*(\mathcal{Y}_t)$ the map on cohomology induced by $\Phi_{\mathcal{P}_t}$. By theorem \ref{thm-equiv-K3}, $\Phi^H_{\mathcal{P}_t}: \widetilde{H}(\mathcal{X}_t,\mathbb{Z})\rightarrow \widetilde{H}(\mathcal{Y}_t,\mathbb{Z})$ is a Hodge isometry. Recall that 
\[
\Phi^H_{\mathcal{P}_t}(\alpha)=q_{t*} (\mathrm{ch}(\mathcal{P}_t)\sqrt{\mathrm{td}(\mathcal{X}_t\times \mathcal{Y}_t)}\cdot p_t^*\alpha).
\]
Since $\mathrm{ch}(\mathcal{P}_t)\sqrt{\mathrm{td}(\mathcal{X}_t\times \mathcal{Y}_t)}$ is a restriction of an algebraic cohomology class on $\mathcal{X}_{\Delta^*}\times_{\Delta^*}\mathcal{Y}_{\Delta^*}$, we have a commutative diagram
\begin{equation}\label{diag-Mukai-1}
\xymatrix{
  H^2(\mathcal{X}_t)\oplus H^0(\mathcal{X}_t)\oplus H^4(\mathcal{X}_t) \ar[r]^{N_X} \ar[d]^{\Phi_{\mathcal{P}_t}^H} &
 H^2(\mathcal{X}_t)\oplus H^0(\mathcal{X}_t)\oplus H^4(\mathcal{X}_t) \ar[d]^{\Phi_{\mathcal{P}_t}^H} \\
 H^2(\mathcal{Y}_t)\oplus H^0(\mathcal{Y}_t)\oplus H^4(\mathcal{Y}_t) \ar[r]^{N_Y}  &
 H^2(\mathcal{Y}_t)\oplus H^0(\mathcal{Y}_t)\oplus H^4(\mathcal{Y}_t) .
}
\end{equation}
So the smallest integer $i$ such that $N_X^i=0$ is equal to that for $N_Y$. We consider the three cases separately.
\begin{enumerate}
  \item[(i)] $N_X=N_Y=0$. By theorem \ref{thm-deg-K3}, $\mathcal{X}_0$ and $\mathcal{Y}_0$ are K3 surfaces.   By proposition \ref{prop-domain}, there is a Hodge isometry between $\widetilde{H}(\mathcal{X}_0,\mathbb{Z})$ and $\widetilde{H}(\mathcal{Y}_0,\mathbb{Z})$. So by theorem \ref{thm-equiv-K3}, $\mathcal{X}_0$ and $\mathcal{Y}_0$ are derived equivalent, so $\rho([X])=\rho([Y])$.
  \item[(ii)] $N_X\neq 0$, $N_Y\neq 0$, $N_X^2=N_Y^2=0$. Then with the notation of theorem \ref{thm-deg-K3}, we have
  \[
  \rho([X])=[V_0]+[V_r]+\sum_{i=1}^{r-1}[V_i]-2r[E]=[V_0]+[V_r]-2[E].
  \]
  Since $e(X)=0$, we have $e(V_0)+e(V_r)-2e(E)=0$, thus $e(V_0)+e(V_r)=0$. Since $V_0$ and $V_r$ are rational surfaces, we have $[V_0]+[V_r]=0$. Therefore
  \[
  \rho([X])=-2[E].
  \]
  It suffices to show $E_X\cong E_Y$. The diagram (\ref{diag-Mukai-1}) induces an isomorphism of Hodge structures
  \[
  N_X(\widetilde{H}(X,\mathbb{Z}))\xrightarrow{\sim} N_Y(\widetilde{H}(Y,\mathbb{Z}))
  \]
  But $N_X(\widetilde{H}(\mathcal{X}_t,\mathbb{Z}))=N_X(H^2(\mathcal{X}_t,\mathbb{Z}))$ and $N_Y(\widetilde{H}(\mathcal{Y}_t,\mathbb{Z}))=N_Y(H^2(\mathcal{Y}_t,\mathbb{Z}))$. By definition of the weight filtration on $LH^2(\mathcal{X}_t)$ and $LH^2(\mathcal{Y}_t)$, $N_X(H^2(\mathcal{X}_t))=W_1 LH^2(\mathcal{X}_t)$ and $N_Y(H^2(\mathcal{Y}_t))=W_1 LH^2(\mathcal{Y}_t)$, as pure Hodge structures. So by proposition \ref{prop-deg-type-ii}, $E_X\cong E_Y$.
  \item[(iii)] $N_X^2\neq 0$, $N_Y^2\neq 0$, $N_X^3=N_Y^3=0$. Since $e(\mathcal{X}_0)=0$ and all the components of $\mathcal{X}_0$ are rational, $\rho([X])=0$. The same holds for $\mathcal{Y}_0$. So we are done.
\end{enumerate}
\pqed

\section{Specialization map for abelian varieties}
In this section we verify (\ref{eq-rho-check-0}) for derived equivalent abelian varieties. In the final result (corollary \ref{cor-ab-dege-1}) we  need to assume that $k$ is an algebraically closed field, because  theorem \ref{thm-pol} need this assumption. However we still state the intermediate statements in a more general setting.
\subsection{Derived equivalences of abelian abelian varieties}
In this subsection we collect some theorems on derived equivalent abelian varieties due to Mukai, Polishchuk and Orlov. Our references are \cite{Muk87b}, \cite{Pol96}, \cite{Orl02}, and also \cite[Chapter 9]{Huy06}.
\begin{theorem}[\cite{Muk87b}]
Let $S$ be a scheme, $p:A\rightarrow S$ an abelian scheme, and $q:A^t\rightarrow S$ its dual abelian scheme:
\[
\xymatrix{
  & A\times_S A^t \ar[dl]_{\pi_A} \ar[dr]^{\pi_{A^t}}  &\\
  A \ar[dr]_{p} && A^t \ar[dl]^{q} \\
  & S&.
}
\]
Denote by $\mathcal{P}$ the Poincaré invertible sheaf on $A\times_S A^t$. Then the Fourier-Mukai functor
\[
\Phi: \mathrm{D}^{\mathrm{b}}(A)\rightarrow \mathrm{D}^{\mathrm{b}}(A^t),\quad \Phi(\mathcal{E})=\mathbf{R}\pi_{A^t*}(\mathbf{L}\pi_A^*(\mathcal{E})\otimes^{\mathbf{L}}\mathcal{P})
\]
is an equivalence.
\end{theorem}

 Let $A$, $B$ be abelian schemes over $S$. Suppose $f:A\times_S A^t\rightarrow B\times_S B^t$ is a homomorphism of abelian varieties. Write $f$ as a matrix
\[
  f=\left(\begin{array}{cc}
  \alpha & \beta \\
  \gamma & \delta
  \end{array}
  \right)
\]
where $\alpha: A\rightarrow B$, $\beta: A^t\rightarrow B$, $\gamma: A\rightarrow B^t$, and $\delta:A^t\rightarrow B^t$. Define a homomorphism $\tilde{f}: B\times_S B^t\rightarrow A\times_S A^t$ by
\[
  \tilde{f}=\left(\begin{array}{cc}
  \delta^t & -\beta^t \\
  -\gamma^t & \alpha^t
  \end{array}
  \right).
\]
\begin{definition}
An isomorphism $f:A\times_S A^t\rightarrow B\times_S B^t$ is called a \emph{symplectic isomorphism} if $f^{-1}=\tilde{f}$.
\end{definition}

\begin{theorem}[\cite{Pol96}]\label{thm-pol}
Let $k$ be  an algebraically closed field,  $A$ and $B$ two abelian varieties over $k$. 
If  there is a sympelctic isomorphism $f:A\times_k A^t\rightarrow B\times_k B^t$, then $A$ and $B$ are derived equivalent.
\end{theorem}

\begin{theorem}[\cite{Orl02}]\label{thm-orl}
Let $k$ be a field, $A$ and $B$ two abelian varieties over $k$.
If $A$ and $B$ are derived equivalent, then there exists a sympelctic isomorphism $f:A\times_k A^t\rightarrow B\times_k B^t$.
\end{theorem}

\subsection{Degeneration and Mumford-Künnemann construction}
From now on, we fix a complete discrete valuation ring $R$, and let $\mathfrak{m}$ be the maximal ideal of $R$, $K$ be the fraction field of $R$, $k$ the residue field of $R$, and denote $S=\Spec(R)$. Denote by $\eta$ and 0 the generic and the closed point of $S$, respectively. 

In this subsection, we recall some notions in the theory of degeneration of abelian varieties. Our references are \cite[chapter 2, 3]{FC90}, \cite[chapter 3, 4]{Lan13}. Then we state a theorem of Künnemann \cite{Kunn98} on the construction of an snc model of an abelian variety over $K$ which admits a split ample degeneration, or called the \emph{Mumford-Künnemann construction} (see also \cite{Mum72}).  

\begin{definition}
Let $A$ be a abelian variety over $K$. A \emph{semistable degeneration} of $A$ over $S$ is a semiabelian scheme $G$ over $S$ with an isomorphism $G_{\eta}\cong A$. By definition, there is an extension 
\[
0\rightarrow T_0\rightarrow G_0\rightarrow A_0\rightarrow 0
\]
where $A_0$ is an abelian variety over $k$, and $T_0$ is a torus over $k$. If $T_0$ is a split torus, $G$ is called a \emph{split degeneration} of $A$.
\end{definition}

\begin{definition}
An \emph{ample degeneration} of $A$ is a pair $(G,\mathscr{L})$ where $G$ is  a semiabelian degeneration of $A$ over $S$ and $\mathscr{L}$ is a cubical invertible sheaf on $G$ such that $\mathscr{L}_\eta$ is ample. 
\end{definition}

In fact the condition implies that $\mathscr{L}$ is relatively ample.

Denote $S_i=\Spec(R/\mathfrak{m}^i)$. For a semiabelian scheme $G$ over $S$, denote $G_{\mathrm{for}}=\lim G\times_S S_i$, and $\mathcal{L}_{\mathrm{for}}$ the corresponding formal completion of $\mathcal{L}$.
For an ample degeneration $(G,\mathcal{L})$, there is the associated \emph{Raynaud extension}
\[
0\rightarrow T\rightarrow \widetilde{G}\xrightarrow{\pi} \widetilde{A}\rightarrow 0,
\]
such that $\widetilde{G}$ is an algebraization of  the formal scheme  $G_{\mathrm{for}}$, $T$ is a torus over $S$, and $\widetilde{A}$ is an abelian scheme over $S$, and there is a cubical ample invertible sheaf $\widetilde{\mathcal{L}}$ which is the algebraization of $\mathcal{L}_{\mathrm{for}}$.

\begin{definition}
A \emph{Split ample degeneration} of $A$ is a triple $(G,\mathcal{L},\mathcal{M})$, where $G$ is a split degeneration of $A$, $(G,\mathcal{L})$ is an ample degeneration, and $\mathcal{M}$ is a cubical ample invertible sheaf on $\widetilde{A}$ such that $\pi^* \mathcal{M}\cong\widetilde{\mathcal{L}}$.
\end{definition}

By the rigidity of tori \cite[X. theorem 3.2]{DG70}, $T_0$ is split implies that $T$ is split. Moreover, the character group of $T$ is a constant abelian sheaf over $S$, and we denote the associated constant group by $X$. There is a notion of dual semiabelian scheme $G^t$ over $S$, and the corresponding torus $T^t$ is also split. We denote the constant character group of $T^t$ by $Y$.

\begin{definition}
Consider the cone $\mathscr{C}=(X_{\mathbb{R}}^{*}\times \mathbb{R}_{>0})\cup \{0\}$. There is a natural action of $Y$ on $\mathscr{C}$ via addition. A \emph{$Y$-admissible polyhedral cone decomposition} of $\mathscr{C}$ is  a (possibly infinite) rational polyhedral cone decomposition $\{\sigma\}_{\alpha\in I}$ of $\mathscr{C}$ such that the collection of the cones $\sigma_\alpha$ is invariant under the action of $Y$ and  there are only finitely many orbits.
\end{definition}

\begin{theorem}\label{thm-Kunn}(\cite[theorem 3.5]{Kunn98})
Let $(G,\mathscr{L},\mathscr{M})$ be a split ample degeneration. Then there is a projective regular model $P$, and an admissible cone decomposition $\{\sigma_{\alpha}\}_{\alpha\in I}$ of $\mathscr{C}$, and we denote by $I_Y^{+}$ the corresponding orbit space with the orbit of the zero cone removed, such that
\begin{enumerate}
  \item[(i)]  the reduced special fiber $(P_0)_{\mathrm{red}}$ is a strict normal crossing divisor on $P$;
  \item[(ii)] $(P_0)_{\mathrm{red}}$ has a natural stratification with strata $G_{\sigma_\alpha}$ for $\alpha\in I_Y^{+}$, where $G_{\sigma_\alpha}$ is a semiabelian scheme fitting into an exact sequence
  \[
0\rightarrow T_{\sigma_{\alpha}}\rightarrow G_{\sigma_{\alpha}}\rightarrow A_0\rightarrow 0,
\]  
where $A_0$ is the abelian part of the Raynaud extension,  and $T_{\sigma_{\alpha}}$ is a split torus;
  \item[(iii)] the closure $P_{\sigma_\alpha}$ of the stratum $G_{\sigma_\alpha}$ is the disjoint union of all $G_{\sigma_\beta}$ such that $\alpha$ is a face of $\beta$, and
  \[
P_{\sigma_{\alpha}}=G_{\sigma_{\alpha}}\times^{T_{\sigma_{\alpha}}}Z_{\sigma_{\alpha}}.
\]
is a contraction product, 
where  $T_{\sigma_{\alpha}}\rightarrow Z_{\sigma_{\alpha}}$ an open torus imbedding into a smooth projective toric  variety.
\end{enumerate}
\end{theorem}

\subsection{Degeneration and derived equivalence}

\begin{proposition}\label{prop-split-ample-1}
 Let $A$ be an abelian varieties over $K$, which has a split  degeneration  over $R$. Then $A$ has a split ample degeneration over $R$.
\end{proposition}
Proof : By the assumption there is  a semi-abelian scheme $G$ over $R$ such that $G_K\cong A$ and $G_0$ fits into an extension 
\[
0\rightarrow T_0\rightarrow G_0\rightarrow A_0\rightarrow 0
\]
such that $T_0$ is a split torus over $k$ and $A_0$ is an abelian variety over $k$. By \cite[I, 2.6]{MB85} and \cite[XI, 1.13]{Ray70} (see also \cite[remark. 3.3.3.9]{Lan13}), there is an ample cubical invertible sheaf $\mathscr{L}$ over $G$. Thus $\mathscr{L}\otimes [-1]^* \mathscr{L}$ is also an ample cubical invertible sheaf over $G$. Let 
\[
0\rightarrow \widetilde{T}\rightarrow \widetilde{G}\rightarrow \widetilde{A}\rightarrow 0
\]
be the corresponding Raynaud extension. Then by \cite[cor. 3.3.3.3, prop. 3.3.3.6]{Lan13} and \cite[XI, 1.11]{Ray70}, the invertible sheaf $\mathscr{L}_{\mathrm{for}}\otimes [-1]^* \mathscr{L}_{\mathrm{for}}$  over $\widetilde{G}_{\mathrm{for}}$ is isomorphic to an ample pullback $\mathscr{M}_{\mathrm{for}}$ over $\widetilde{A}_{\mathrm{for}}$ which is algebraizable. This provides a split ample degeneration of $A$.\pqed
\begin{lemma}\label{lem-extension-1}
Let 
\[
0\rightarrow T\rightarrow G\rightarrow A\rightarrow 0
\]
be an extension of an abelian variety $A$ by a split torus, over a field $k$, and $T\hookrightarrow Z$ be an open torus embedding of $T$ into a smooth complete toric variety $Z$.
 Then in $K_0(\mathrm{Var}_k)$ one has
 \[
 [G\times^T Z]=[Z]\cdot [A].
 \]
\end{lemma}
Proof : By the assumption, $G$ is an fppf $T$-torsor over $A$. Since $T$ is split, $G$ is a product of fppf $\mathbb{G}_m$-torsors over $A$, thus it is also a product of  Zariski $\mathbb{G}_m$-torsors, by  Hilbert theorem 90. So there is a locally closed stratification $\{U_\alpha\}$ of $A$ such that $(G\times^T Z)|_{U_\alpha}$ is isomorphic to $Z\times U_\alpha$, hence the conclusion.\pqed

\begin{theorem}\label{thm-ab-dege-2}
Let $(R,\mathfrak{m})$ be a complete discrete valuation ring, $K$ the fraction field of $R$, and $k$ the residue field of $R$. Let $A$ and $B$ be two abelian varieties over $K$, which are derived equivalent.  Then the following holds.
\begin{enumerate}
  \item[(i)] $A$ has semistable reduction if and only if $B$ has semistable reduction.
  \item[(ii)] $A$ has a split  degeneration if and only if 
  $B$ has a split  degeneration. 
  \item[(iii)] In case of (i), denote the abelian part of the special fiber of the semistable reduction of $A$ (resp., $B$) by $A_0$ (resp., $B_0$). Then there is a symplectic isomorphism $A_0\times_k A_0^t\xrightarrow{\sim} B_0\times_k B_0^t$.
\end{enumerate}
\end{theorem}
Proof: (i) By theorem \ref{thm-orl} there is a symplectic isomorphism $A\times_K A^t\cong B\times_K B^t$. Then by \cite[proposition 2.10]{MP17}, $A$ and $B$ are isogenous. Thus the conclusion (i) follows from \cite[\S 7.3, corollary 3]{BLR90}.

(ii) Suppose $A$ has a split ample degeneration over $R$. By \cite[\S 2.2]{FC90}, $A^t$ has a split ample degeneration over $R$. Let $\mathscr{A}$ (resp., $\mathscr{A}'$) be the Néron model of $A$ (resp., of $A^t$) over $R$. By the functoriality of Néron models, $\mathscr{A}\times_R \mathscr{A}'$ is the Néron model of $A\times_K A^t$. By theorem \ref{thm-orl}, $A\times_K A^t$ is isomorphic to $B\times_K B^t$. Let $\mathscr{B}$ (resp., $\mathscr{B}'$) be the Néron model of $B$ (resp., of $B^t$) over $R$. Thus $\mathscr{A}\times_R \mathscr{A}'\cong \mathscr{B}\times_R \mathscr{B}'$, so their special fibers have isomorphic  identity components, i.e. $(\mathscr{A}\times_R \mathscr{A}')_k^{\circ}\cong (\mathscr{B}\times_R \mathscr{B}')_k^{\circ}$. By (i), $A$, $A^t$, $B$, $B^t$ all have semistable reductions over $R$. Thus $(\mathscr{A}_k)^{\circ}$, $(\mathscr{A}')^{\circ}$, $(\mathscr{B})^{\circ}$ and $(\mathscr{B}')^{\circ}$ are all semi-abelian varieties over $k$, hence are geometrically connected. Thus by \cite[4.5.8]{EGAIV} $(\mathscr{A}_k)^{\circ}\times_k (\mathscr{A}_k')^{\circ}$ and $(\mathscr{B}_k)^{\circ}\times_k (\mathscr{B}_k')^{\circ}$ are connected and thus are isomorphic to $(\mathscr{A}\times_R \mathscr{A}')_k^{\circ}$. Let $T$ (resp. $T'$) be the torus part of $\mathscr{B}_k$ (resp. $\mathscr{B}_k'$). Then $T\times_k T'$ is a split torus. Consider the character group $\underline{X}(T)$ (resp. $\underline{X}(T')$) of $T$ (resp. $T'$), which are étale sheaves of torsion free abelian groups of finite type. The product $\underline{X}(T)\times \underline{X}(T')$ is the character group of $T\times_k T'$, and is therefore a constant sheaf by the splitness of $T\times_k T'$. Considering the action of $\mathrm{Gal}(k^s/k)$ on $\underline{X}(T)(k^s)$ and $ \underline{X}(T')(k^s)$, one sees that both $\underline{X}(T)$ and $\underline{X}(T')$ are constant sheaves over $k_{\mathrm{\acute{e}t}}$, and therefore $T$ and $T'$ are split tori over $k$.

(iii) By theorem \ref{thm-orl} there is an isomorphism $f:A\times_K A^t\xrightarrow{\sim}B\times_K B^t$ of the form
\[
  f=\left(\begin{array}{cc}
  \alpha & \beta \\
  \gamma & \delta
  \end{array}
  \right)
\]
such that 
\[
\left(\begin{array}{cc}
  \delta^t & -\beta^t \\
  -\gamma^t & \alpha^t
  \end{array}
  \right)\cdot \left(\begin{array}{cc}
  \alpha & \beta \\
  \gamma & \delta
  \end{array}
  \right)=\id.
\]
By the functoriality of Néron models the  isomorphism $f$ extends to an isomorphism $F: \mathscr{A}\times_R \mathscr{A}'\xrightarrow{\sim} \mathscr{B}\times_R \mathscr{B}'$ of the form
\[
  F=\left(\begin{array}{cc}
  \widetilde{\alpha} & \widetilde{\beta} \\
  \widetilde{\gamma} & \widetilde{\delta}
  \end{array}
  \right)
\]
such that 
\[
\left(\begin{array}{cc}
  \widetilde{\delta}^t & -\widetilde{\beta}^t \\
  -\widetilde{\gamma}^t & \widetilde{\alpha}^t
  \end{array}
  \right)\cdot \left(\begin{array}{cc}
  \widetilde{\alpha} & \widetilde{\beta} \\
  \widetilde{\gamma} & \widetilde{\delta}
  \end{array}
  \right)=\id.
\]
Considering the special fibers and using the proof of (ii), one obtains a symplectic isomorphism between the abelian parts of $\mathscr{A}_k^{\circ}\times_k (\mathscr{A}_k')^{\circ}$ and $\mathscr{B}_k^{\circ}\times_k (\mathscr{B}_k')^{\circ}$.
\pqed

\begin{proposition}\label{prop-ab-dege-2}
Let $A$ and $B$ be two abelian varieties over $K$, which are derived equivalent, and suppose that $A$ has a split degeneration over $R$. Then $A$ and $B$ have snc models $P$ and $Q$ over $R$, respectively, such that either $P_0$ and $Q_0$ are symplectically isomorphic abelian varieties over $k$, or $[P_0]=[Q_0]=0$ in $K_0(\mathrm{sGT}_k)$.
\end{proposition}
Proof : By theorem \ref{thm-ab-dege-2} (ii) and proposition \ref{prop-split-ample-1}, both $A$ and  $B$ has split ample degeneration over $R$. By theorem \ref{thm-ab-dege-2} (iii), if $A$ has good reduction over $R$, then so does $B$, and $A_0$ and $B_0$ are symplectic isomorphic. If $A$ does not have a good reduction over $R$, then by theorem \ref{thm-Kunn} and lemma \ref{lem-extension-1}, 
\[
[P_0]=\sum_{a\in I_Y^+} (-1)^{j_\alpha-1}j_{\alpha}
[A_0]\times [Z_{\sigma_\alpha}]
\]
in $K_0(\mathrm{Var}_k)$, where 
\[
j_\alpha=\dim A+1-\dim A_0-\dim Z_{\sigma_\alpha}=\dim_{\mathbb{R}}\mathscr{C}-\dim_{\mathbb{R}}\sigma_\alpha.
\]
Since each face of $\sigma_\alpha$ appears in the above sum, a simple manipulation shows that
\[
\sum_{a\in I_Y^+} (-1)^{j_\alpha-1}j_{\alpha} [Z_{\sigma_\alpha}]
\]
is equal to a linear combination of split tori in $K_0(\mathrm{Var}_k)$, so
\[
[P_0]\equiv 0\mod (\mathbb{L}-1).
\]
Then by corollary \ref{cor-mu} and the definition (\ref{eq-def-rho-sgt-1}) of $\rho$ , $[P_0]=0$ in $K_0(\mathrm{sGT}_k)$. By theorem \ref{thm-ab-dege-2}, $B$ also has a split ample but not good degeneration over $R$, thus one has $[Q_0]=0$ in $K_0(\mathrm{sGT}_k)$, too. So we are done.\pqed

\begin{corollary}\label{cor-ab-dege-1}
Let $(R,\mathfrak{m})$ be a complete discrete valuation ring, $K$ the fraction field of $R$, and $k$ the residue field of $R$. Suppose that $k$ is algebraically closed of characteristic 0. Let $A$ and $B$ be two abelian varieties over $K$, which are derived equivalent. Suppose $A$ has a semistable reduction over $R$. Then $\rho([A])=\rho([B])$.
\end{corollary}
Proof : By theorem \ref{thm-ab-dege-2} (i), both $A$ and $B$ semistable reductions over $R$, which are automatically split degenerations because $k$ is algebraically closed. Applying proposition \ref{prop-ab-dege-2} and theorem \ref{thm-pol} we obtain the conclusion. \pqed

\section{Open problems}
\begin{enumerate}
  \item Although  our (conjectural) definition of $\rho_{\mathrm{sgt}}$ does not assume the existence of semistable degeneration over $R$, in the above verifications we need to assume this to apply the results for the degeneration of these varieties. It is natural to make the following conjecture. Theorem \ref{thm-ab-dege-2} provides an example for it.

\begin{conjecture}
Let $R$ be a DVR, $K$ its fraction field. Let $X$ and $Y$ be derived equivalent smooth projective varieties over $K$. Then $X$ has semistable  degeneration (resp., good reduction) over $R$ if and only if $Y$ has semistable degeneration (resp., good reduction) over $R$.
\end{conjecture}
This suggests to take into consideration the Galois action on the derived categories, and ask whether there is a Néron-Ogg-Shafarevich-Grothendieck type criterion for the types of degenerations.
\item Does there exist a smooth projective variety $X$ over $k$ such that $[X]=m[\Spec(k)]$ in $K_0(\mathrm{sGT})$ but $\mathrm{D}^{\mathrm{b}}(X)$ does not have a full exceptional collection? If there are such varieties, are their quantum cohomology semisimple? The limit fibers of a family of varieties with full exceptional collections are candidates for this.
\end{enumerate}

\textsc{School of Mathematics, Sun Yat-sen University, Guangzhou 510275, P.R. China}

\emph{Email address}: huxw06@gmail.com

\end{document}